\documentclass[a4paper]{article}
\usepackage{amsmath,bm,stmaryrd}
\usepackage{amssymb}
\usepackage{graphicx}
\usepackage{epstopdf}
\usepackage{multirow}
\usepackage{calc}
\usepackage[margin=1.in]{geometry}
\usepackage{amsmath, amsthm, amssymb}
\usepackage[usenames]{color}
\usepackage{indentfirst}
\usepackage[shortlabels]{enumitem}
\usepackage{float}
\usepackage[mathscr]{euscript}
\usepackage{mathtools}

\usepackage[figuresright]{rotating}
\usepackage[ruled,vlined]{algorithm2e}
\usepackage{soul}
\usepackage[title]{appendix}

\usepackage[colorlinks,
            linkcolor=blue,
            anchorcolor=blue,
            citecolor=red
            ]{hyperref}
\usepackage{booktabs}

\newcommand\keywordsname{Key words}
\newcommand\AMSname{AMS subject classifications}

\newenvironment{@abssec}[1]
{\if@twocolumn
\section*{#1}%
\else
\vspace{.05in}\footnotesize
\parindent .2in
{\upshape\bfseries #1. }\ignorespaces
\fi}

{\if@twocolumn\else\par\vspace{.1in}\fi}

\newenvironment{keywords}{\begin{@abssec}{\keywordsname}}{\end{@abssec}}

\providecommand{\Div}{\operatorname{div}}          

\providecommand{\Dim}{\operatorname{dim}}            
\providecommand{\dim}{\Dim}

\providecommand*{\Dist}[2]{\operatorname{dist}({#1};{#2})}   
\providecommand*{\Dist}[2]{\Dist{#1}{#2}}

\newcommand{\be}{\begin{eqnarray}}
\newcommand{\ee}{\end{eqnarray}}

\newcommand{\ben}{\begin{eqnarray*}}
\newcommand{\een}{\end{eqnarray*}}

\title{Robust Block Preconditioning for 3D nonlinear steady-state radiation transport equations}
\author{Yunpan Ma \thanks{School of Mathematics and Statistics, Zhengzhou University, Zhengzhou, 450001, Henan, China.}
\and
Lingxiao Li
\thanks{School of Mathematics and Statistics, Henan University, Kaifeng, 475004, Henan, China.}
\and
Changhui Yao
\thanks{School of Mathematics and Statistics, Zhengzhou University, Zhengzhou, 450001, Henan, China}}
\begin{document}
\date{}
\maketitle

\begin{abstract}
  In this work, based on the discrete ordinate method,
  we propose a robust block preconditioning strategy for the 3D nonlinear steady-state
  radiation transport equation with heat diffusion term. The presence of the diffusive term of the temperature equation prevents its elimination into a single equation for the radiation intensity. To overcome this difficulty,
  all physical variables are assembled into a single monolithic linear system. The heat flux and temperature are
  treated as independent variables in a mixed $H(\mathrm{div})$-conforming finite element formulation.
  The equation for radiation intensity
  is discretised by a discontinuous Galerkin method with upwind flux, where a vectorial finite element space is used to couples the radiation intensity in different directions within each element.
  We then construct a Newton-Krylov iterative solver to solve the nonlinear equations,
  for which the core part is efficient preconditioning.
  To accelerate the convergence of Krylov's method, three block preconditioners are constructed,
  corresponding to different levels of approximation of the coupling between the temperature and radiation intensity. $P_{\mathrm{Schur}}$ retains the full coupling. $P_{\mathrm{Split}}$ drops the conductive contribution to the radiation block. $P_{\mathrm{BJ}}$ neglects the radiation-to-temperature coupling, retaining only the temperature-to-radiation coupling. Numerical experiments demonstrate the mesh independence and robustness of the proposed preconditioners.
\end{abstract}


\section{Introduction}

Modeling and simulating heat radiation transport is crucial in a wide range of high-temperature applications, including inertial confinement fusion (ICF)~\cite{Haines2024}, magnetic confinement fusion (MCF)~\cite{Sciortino2021}, and astrophysical phenomena~\cite{Wunsch2024}. The key physical processes governing energy transfer in these applications are radiation transport and thermal conduction, whose strong coupling poses significant challenges for numerical simulations~\cite{Modest2013, Howell2020}. In this work, we consider the three-dimensional steady-state system described by the radiation transport equation (RTE) and the energy conservation equation
\begin{align}
  \Omega\cdot\nabla\Phi + \sigma_a\Phi - \sigma_a B(T) &= f \quad \text{in } D\times \mathbb{S}^{2}, \label{eq:intro_rte}\\
  -\nabla\cdot(k\nabla T) - \sigma_a\int_{\mathbb{S}^{2}}\Phi\,\mathrm{d}\Omega + 4\pi\sigma_a B(T) &= 0 \quad \text{in } D, \label{eq:intro_energy}
\end{align}
where $\Phi(\mathbf{x},\Omega)$ is the radiative intensity, $T(\mathbf{x})$ the temperature, and $B(T)=\sigma_B T^4/\pi$ the Planck function. The coefficients $\sigma_B$, $\sigma_a$, and $k$ denote the Stefan-Boltzmann constant, the absorption coefficient, and the thermal conductivity, respectively.
The RTE is a five-dimensional integro-differential equation. Its high dimensionality arises from the dependence of the radiative intensity $\Phi$ on three spatial and two angular coordinates. {Moreover, not that the presence of the heat diffusive term $-\nabla\cdot(k\nabla T)$ in the energy equation couples the temperature directly to the radiation field,
preventing the simple use of single equation for radiation intensity to implement the traditional sweeping algorithm.}
After discretisation, this coupling yields a nonsymmetric and indefinite linear system, for which iterative solvers without preconditioning converge slowly or even fail. The main contribution of this work is the development of robust preconditioners for linear systems arising from numerical discretisation of this coupled system.

The numerical solution of the RTE has been extensively studied over the past decades. Solution strategies are broadly classified into stochastic and deterministic approaches. Stochastic methods, primarily Monte Carlo techniques~\cite{Shi2023a, Shi2023b}, simulate photon transport via random particle histories. They offer geometric flexibility and accurate physics modelling. {Unfortunately, the slow statistical convergence renders them prohibitively expensive for optically thick or large-scale problems.} Deterministic methods~\cite{Hang2014} are typically represented by the discrete ordinates method ($S_N$ method)~\cite{CarlsonLathrop1968, Coelho2014} and spherical harmonics ($P_N$ method) approximations~\cite{Pomraning1979}. These discretisations are commonly solved using source iteration (SI)~\cite{LewisMiller1984, AdamsLarsen2002}, an operator-split strategy that alternates between transport sweeps and updates of scattering and emission sources. The convergence of source iteration was first analysed by Reed~\cite{Reed1971}. Fourier analysis~\cite{Larsen1988} later showed that its spectral radius is governed by the scattering ratio and is independent of the spatial mesh. Strong radiation-matter coupling and reflective boundaries can further degrade convergence in practice, motivating the development of acceleration techniques such as diffusion synthetic acceleration (DSA)~\cite{Alcouffe1977}, grey transport acceleration (GTA)~\cite{Larsen1988}, and transport synthetic acceleration (TSA)~\cite{Ramone1997}. More recently, Tang et al.~\cite{Tang2021} developed a semi-implicit asymptotic preserving scheme for the gray RTE that avoids nonlinear iterations; a frequency-dependent extension was later proposed by Zhang et al.~\cite{Zhang2023}. Beyond these acceleration strategies, alternative approaches have also been explored, including angular-spatial multigrid methods~\cite{GaoZhao2009, GaoZhao2013}, asymptotic-preserving discretisations~\cite{GuermondKanschat2010}, unified gas kinetic schemes~\cite{Tan2019}, and flux-corrected transport techniques~\cite{HanselRagusa2018}.

The effectiveness of DSA, GTA, and TSA relies on the operator-split structure of SI, in which the temperature appears only as a source term updated after each transport sweep. {Unfortunately, the diffusive term $-\nabla\cdot(k\nabla T)$ directly couples the temperature and radiation fields, rendering independent updates inappropriate.}
Decoupled alternating updates converge slowly because the thermal diffusion time scale is much longer than that of radiative transport, thus precluding an effective operator-split strategy, and such treatment may imposes a prohibitive CFL stability restriction on the time step size.

To overcome these difficulties, Krylov subspace methods provide a flexible framework for solving the linear systems arising from discrete ordinates discretisations~\cite{Charest2012, GodoyLiu2012}. Existing implementations follow two distinct strategies.
\begin{itemize}
\item Use Krylov methods as an outer accelerator for the standard transport sweep, with synthetic acceleration as a preconditioner.
\item Assemble all discrete directions into a monolithic matrix and solve with preconditioned Krylov methods.
\end{itemize}
These two strategies offer complementary strengths. Sweep-based methods are computationally efficient and have been extensively studied in the literature. Synthetic acceleration techniques, such as linear multifrequency-grey (LMFG) acceleration~\cite{Till2018} and reduced-order model (ROM)-enhanced acceleration~\cite{Peng2025FGMRES}, have been used as effective preconditioners for Krylov methods. DSA~\cite{Warsa2004} and TSA~\cite{Drumm2004, DrummFan2017} have also been applied both as stand-alone accelerators for source iteration and as preconditioners for Krylov methods. Monolithic strategies~\cite{Badri2018, Badri2019, Li2021, Jolivet2021}, in contrast, eliminate the need for a topological sweep order, thereby circumventing the scheduling difficulties associated with sweeping-based methods on complex unstructured grids (or even curved meshes)~\cite{Plimpton2005}. This makes them particularly attractive for complex geometries and heterogeneous media. Moreover, the linear system can be distributed without complex sweep scheduling, and the thermal conduction term in the present coupled system makes the monolithic approach a natural choice by avoiding the convergence difficulties inherent in operator splitting.

The spatial discretisation of the RTE is typically performed using finite volume methods (FVM) or finite element methods (FEM). FVM remains a common choice in engineering radiative transfer due to its local conservation property: the exact balance of radiative energy at the control volume level and its ease of implementation on structured grids~\cite{Coelho2014}. Unfortunately, its accuracy on unstructured meshes is typically limited to first or second order. Continuous FEM offers greater geometric flexibility than FVM, but requires additional stabilisation (e.g., SUPG~\cite{Badri2018, Badri2019, LeHardy2016}) to suppress spurious oscillations when applied to the hyperbolic RTE. Discontinuous Galerkin (DG) methods are a natural choice for unstructured meshes~\cite{CockburnShu1998}. They offer high-order accuracy and capture the hyperbolic character of radiation transport via upwind numerical fluxes~\cite{GaoZhao2013, Houston2002}. Moreover, DG methods with upwind fluxes are asymptotically preserving in the diffusive limit, ensuring accurate solutions across the full range of optical thicknesses~\cite{GuermondKanschat2010}, with hp-DG methods providing robust stabilization for first-order hyperbolic problems~\cite{Houston2000, Houston2002}. In the coupled conduction-radiation system, material properties such as absorption, scattering, and thermal conductivity can vary sharply across the domain, leading to steep solution gradients that DG methods resolve stably. Positivity preservation is essential for physically meaningful solutions, and high-order DG methods with positivity-preserving limiters guarantee non-negative radiative intensity~\cite{Yuan2016}. Adaptive mesh refinement (AMR) can further improve computational efficiency by locally concentrating degrees of freedom in regions with steep gradients or sharp layers~\cite{Zhang2020}. To achieve local energy conservation at the element level, we treat the heat flux $\mathbf{q}$ as an independent variable in a mixed finite element formulation~\cite{Brezzi1974, RaviartThomas1977, BrezziFortin1991}, discretised in $H(\Div)$-conforming spaces. This ensures that the divergence of $\mathbf{q}$ balances the radiation source term locally, rather than being recovered through post-processing from the temperature gradient. The resulting discrete system inherits a saddle-point block structure~\cite{Pazner2024}, for which we exploit in the design of our preconditioner.

In this paper, we develop a family of preconditioned Krylov methods for the monolithic linear system~\cite{Badri2018, Badri2019, Li2021} arising from the coupled nonlinear radiation transport system \eqref{eq:intro_rte}
and \eqref{eq:intro_energy}. To handle the $T^4$ nonlinearity of the Planck function, we linearise the coupled system via Newton's method~\cite{Mousseau2000}. The angular variable is discretised using $S_N$ method. The spatial discretisation treats the heat flux and temperature as independent variables for the conduction equation and employs DG with upwind fluxes for the radiation intensity. Based on this discretisation, we assemble all discrete directions and spatial degrees of freedom into a single monolithic linear system. The radiation transport part is assembled using a vectorial finite element
method that couples angular directions within the same spatial element~\cite{Badri2018Vectorial,Li2021}.

The main contribution of this article is the development of robust block preconditioners for the monolithic system.
We apply a approximate Schur complement to eliminate the flux-temperature-coupled block~\cite{BrownWoodward2001, Yue2021, Feng2012, Pazner2024}, with the inverse of the $H(\Div)$ mass matrix approximated by its diagonal. Based on this reduced system, we construct a hierarchy of three preconditioners that represent different approximations to the coupling between the temperature and radiation blocks.
$P_{\mathrm{Schur}}$ approximates the full coupling structure most accurately among the three and provides consistent performance across all regimes. $P_{\mathrm{Split}}$ drops the conductive contribution from the radiation block while keeping the radiation-temperature coupling, yielding a simplified approximation that remains effective in the radiation-dominated regime. $P_{\mathrm{BJ}}$ further neglects the radiation-to-temperature coupling, while retaining the temperature-to-radiation coupling, reducing to the simplest preconditioner in the hierarchy, which is effective when the coupling is weak-i.e. in the conduction-dominated regime.
Their performance is evaluated across conduction-dominated, radiation-dominated, and strong-scattering regimes.

This paper is organized as follows. Section~2 describes a dimensionless coupled nonlinear heat
radiation transport equations and its Newton linearisation.
Section~3 presents the discrete ordinates method and the mixed DG finite element discretisation.
Three block preconditioners are constructed to develop a preconditioned FGMRES algorithm for solving the algebraic systems in Section~4. Numerical experiments demonstrating convergence and preconditioner performance are reported in Section~5. Concluding remarks are given in Section~6.

\section{Model equations and Newton linearisation}

In this section we present the nondimensional governing equations and their Newton linearisation, which form the basis for the numerical discretisation in Section~\ref{sec:discretisation}.

\subsection{Nondimensional model}

To ensure local energy conservation, we introduce the heat flux $\mathbf{q}$ as an auxiliary variable. Following the
derivation detailed in~\ref{sec:appendices}, the dimensionless steady-state coupled nonlinear transport equation reads
as follows:
\begin{align}
\kappa^{-1}\mathbf{q} + \nabla T &= \mathbf{0} \quad \text{in } D, \label{eq:model1}\\
\operatorname{div}\mathbf{q} - \beta\sigma\left(\frac{1}{4\pi}\int_{\mathbb{S}^{2}}\Phi\,\mathrm{d}\Omega - T^4\right) &= 0 \quad \text{in } D, \label{eq:model2}\\
\Omega\cdot\nabla\Phi + \sigma\Phi - \sigma T^4 &= f \quad \text{in } D\times \mathbb{S}^{2}, \label{eq:model3}
\end{align}
Here $\Phi(\mathbf{x},\Omega)$ denotes the dimensionless radiative intensity, $T(\mathbf{x})$ the dimensionless temperature, and $\mathbf{q}(\mathbf{x})$ the dimensionless heat flux. The dimensionless parameters $\sigma$, $\beta$, and $\kappa$ represent the optical thickness, the radiation-conduction coupling strength, and the dimensionless thermal conductivity, respectively.

The introduction of $\mathbf{q}$ as an independent variable enables a mixed finite element formulation in $H(\mathrm{div})$-conforming spaces, which ensures local energy conservation at the element level. This saddle-point structure will also be exploited in the design of the preconditioner.

The boundary conditions are
\begin{equation}
T = T_D \quad \text{on } \partial D, \qquad
\Phi = \Phi^{\mathrm{in}} \quad \text{on } \Gamma_{\mathrm{in}},
\end{equation}
with $\Gamma_{\mathrm{in}} = \{(\mathbf{x},\Omega)\in\partial D\times\mathbb{S}^{2}: \Omega\cdot\mathbf{n} < 0\}$, and $\mathbf{n}$ the outward unit normal.

\subsection{Newton linearisation}

Given the current iterate \((\mathbf{q}^k,T^k,\Phi^k)\), we employ Newton's method to linearise the system (\ref{eq:model1})-(\ref{eq:model3}) and define the next iterate as
\begin{equation}
  \mathbf{q}^{k+1} = \mathbf{q}^k + \delta\mathbf{q},\qquad
  T^{k+1} = T^k + \delta T,\qquad
  \Phi^{k+1} = \Phi^k + \delta\Phi,
\end{equation}
where the increments $(\delta\mathbf{q},\delta T,\delta\Phi)$ are assumed to be small. Inserting this ansatz into \eqref{eq:model1}-\eqref{eq:model3} and neglecting second-order terms gives the linearised equations for the increments:
\begin{align}
  \kappa^{-1}\delta\mathbf{q} + \nabla\delta T &= R_{\mathbf{q}}^k\quad \text{in } D, \label{eq:newton1}\\
  \operatorname{div}\delta\mathbf{q} - \beta\sigma\left(\frac{1}{4\pi}\int_{\mathbb{S}^{2}}\delta\Phi\,\mathrm{d}\Omega - 4(T^k)^3\delta T\right)
  &= R_T^k\quad \text{in } D, \label{eq:newton2}\\
  \Omega\cdot\nabla\delta\Phi + \sigma\delta\Phi - 4\sigma(T^k)^3\delta T
  &= R_\Phi^k\quad \text{in } D\times \mathbb{S}^{2}, \label{eq:newton3}
\end{align}
where the right-hand side residuals are defined as
\begin{align}
  R_{\mathbf{q}}^k &= -\bigl(\kappa^{-1}\mathbf{q}^k + \nabla T^k\bigr), \\
  R_T^k &= -\operatorname{div}\mathbf{q}^k + \beta\sigma\left(\frac{1}{4\pi}\int_{\mathbb{S}^{2}}\Phi^k\,\mathrm{d}\Omega - (T^k)^4\right), \\
  R_\Phi^k &= f - \bigl(\Omega\cdot\nabla\Phi^k + \sigma\Phi^k - \sigma(T^k)^4\bigr).
\end{align}

Assuming the current iterate satisfies the boundary conditions, the increments inherit homogeneous boundary conditions
\begin{equation}\label{eq:boundary}
  \delta T = 0 \quad \text{on } \partial D, \qquad
  \delta\Phi = 0 \quad \text{on } \Gamma_{\mathrm{in}}.
\end{equation}

Define the combined residual norm $\mathcal{R}^k$ by
\[
  \mathcal{R}^k := \sqrt{\|R_{\mathbf{q}}^k\|_{0,D}^2 + \|R_T^k\|_{0,D}^2 + \|R_\Phi^k\|_{0,D}^2}.
\]

The overall Newton method proceeds as follows:
\begin{enumerate}
  \item Choose an initial guess $(\mathbf{q}^0,T^0,\Phi^0)$ and set $k=0$. Typical choices include $\mathbf{q}^0=\mathbf{0}$, $T^0=0$, $\Phi^0=0$. Compute $\mathcal{R}^0$.
  \item Solve \eqref{eq:newton1}-\eqref{eq:newton3} for the increments.
  \item Update $(\mathbf{q}^{k+1},T^{k+1},\Phi^{k+1}) = (\mathbf{q}^k,T^k,\Phi^k) + (\delta\mathbf{q},\delta T,\delta\Phi)$.
  \item Compute $\mathcal{R}^{k+1}$. If $\mathcal{R}^{k+1} / \mathcal{R}^0 < \epsilon$, stop.
  \item Set $k \leftarrow k+1$ and return to step~2.
\end{enumerate}

Once discretised, each Newton step requires the solution of a large sparse linear system.
The structure of this system will be exploited in the devise of the block preconditioners in Section~4.
The discretisation via the discrete ordinates method and discontinuous Galerkin finite elements is described in the next section.

\section{Angle and finite element discretisation}
\label{sec:discretisation}
In this section, we discretise the linearised integro-differential system \eqref{eq:newton1}-\eqref{eq:newton3}. The discrete ordinates method~\cite{CarlsonLathrop1968, Coelho2014} is used for the angular variable, while the spatial discretisation employs mixed and DG finite elements.

Let $D\subset\mathbb{R}^3$ be a bounded Lipschitz domain with boundary $\partial D$ and outward unit normal $\mathbf{n}$. Let $L^2(D)$ be the usual Hilbert space of square integrable functions, we shall use the following function spaces:
\begin{align*}
H^1(D) &= \left\{ v \in L^2(D) : \nabla v \in [L^2(D)]^3 \right\},\\
\mathbf{H}(\mathrm{div},D) &= \left\{ \mathbf{v} \in [L^2(D)]^3 : \nabla\cdot\mathbf{v} \in L^2(D) \right\}.
\end{align*}

Let $(\cdot,\cdot)_D$ and $\langle\cdot,\cdot\rangle_{\partial D}$ denote the $L^2$ inner products on $D$ and $\partial D$, respectively:
\[
(\mathbf{u},\mathbf{v})_D = \int_D \mathbf{u}\!\cdot\!\mathbf{v}\,\mathrm{d}\mathbf{x},\quad
(u,v)_D = \int_D u v\,\mathrm{d}\mathbf{x},\quad
\langle u,v\rangle_{\partial D} = \int_{\partial D} u v\,\mathrm{d}s,\quad
\|u\|_{0} = (u,u)_D^{1/2}.
\]
For each element $K\in\mathcal{T}_h$, the element-wise inner product $(u,v)_K$ is defined analogously.

For the upwind numerical flux used in the DG discretisation, we introduce the mesh-dependent face inner product
\[
\langle u,v\rangle_{\partial K,\Omega_m} = \int_{\partial K} (\Omega_m\!\cdot\!\mathbf{n})\,u\,v\,\mathrm{d}s,
\]
where $\partial K$ is the boundary of an element $K\in\mathcal{T}_h$ and $\mathbf{n}$ is its outward unit normal. Let $\|M\|$ denote the Euclidean norm of the matrix $M$. These conventions will be used throughout the remainder of the paper.

\subsection{Discrete ordinates method}

Let $\{\Omega_m\}_{m=1}^{N_a}\subset S^2$ be a set of discrete ordinates~\cite{CarlsonLathrop1968, Coelho2014} and $\{w_m\}_{m=1}^{N_a}$ the associated quadrature weights such that
\[
\frac{1}{4\pi}\oint_{S^2} \delta\Phi(\mathbf{x},\Omega)\,\mathrm{d}\Omega
\;\approx\; \sum_{m=1}^{N_a} w_m\,\delta\Phi_m(\mathbf{x}),\qquad
\sum_{m=1}^{N_a} w_m = 1,
\]
with $\delta\Phi_m(\mathbf{x}) = \delta\Phi(\mathbf{x},\Omega_m)$. One can also discretise the angular variable using a quasi-uniform surface mesh of the unit sphere, as in~\cite{Badri2018}.
Evaluating \eqref{eq:newton3} at each $\Omega_m$ gives the semi-discrete system
\begin{align}
  \kappa^{-1}\delta\mathbf{q} + \nabla\delta T &= r_1^{(k)} \quad \text{in } D, \label{eq:sd1}\\
  \operatorname{div}\delta\mathbf{q} - \beta\sigma\Bigl(\sum_{m=1}^{N_a} w_m \delta\Phi_m - 4(T^k)^3\delta T\Bigr) &= r_2^{(k)} \quad \text{in } D, \label{eq:sd2}\\
  \Omega_m\!\cdot\!\nabla\delta\Phi_m + \sigma\delta\Phi_m - 4\sigma(T^k)^3\delta T &= r_{3,m}^{(k)} \quad \text{in } D, \label{eq:sd3}
\end{align}
where $m=1,\dots,N_a$ and the right-hand sides are given by
\begin{align}
  r_1^{(k)} &= -\bigl(\kappa^{-1}\mathbf{q}^k + \nabla T^k\bigr),\\
  r_2^{(k)} &= - \operatorname{div}\mathbf{q}^k + \beta\sigma\Bigl(\sum_{m=1}^{N_a} w_m \Phi_m^{(k)} - (T^k)^4\Bigr),\\
  r_{3,m}^{(k)} &= f_m -\bigl(\Omega_m\cdot\nabla\Phi^k_m + \sigma\Phi^k_m - \sigma(T^k)^4\bigr).
\end{align}
The right-hand sides $r_1^{(k)}$, $r_2^{(k)}$, and $r_{3,m}^{(k)}$ are the residuals of the three governing equations evaluated at the current iterate.

\subsection{Continuous variational formulation}

This subsection presents the continuous variational formulation of the angular semi-discrete problem,
which serves as the foundation for the mixed and DG discretisations in the following.

The continuous mixed variational formulation reads: given the current Newton iterate $(\mathbf{q}^k, T^k, \Phi_m^k)$, find the increments $(\delta\mathbf{q}, \delta T, \delta\Phi_m) \in \mathbf{H}(\mathrm{div},D) \times L^2(D) \times H^1(D)$ such that
\begin{align}
(\kappa^{-1}\delta\mathbf{q},\mathbf{v})_D - (\delta T,\nabla\cdot\mathbf{v})_D
&= R_1(\mathbf{v}) && \forall \mathbf{v}\in \mathbf{H}(\mathrm{div},D), \label{eq:cont1} \\
(\nabla\cdot\delta\mathbf{q},\varphi)_D
+ 4\beta\sigma((T^k)^3\delta T,\varphi)_D
- \beta\sigma\sum_{m=1}^{N_a} w_m (\delta\Phi_m,\varphi)_D
&= R_2(\varphi) && \forall \varphi\in L^2(D), \label{eq:cont2} \\
(\Omega_m\cdot\nabla\delta\Phi_m,\psi)_D + \sigma(\delta\Phi_m,\psi)_D - 4\sigma((T^k)^3\delta T,\psi)_D
&= R_{3,m}(\psi) && \forall \psi\in L^2(D), \label{eq:cont3}
\end{align}
where $m=1,\dots,N_a$, and the residual functionals are defined by
\begin{align}
R_1(\mathbf{v}) &= -(\kappa^{-1}\mathbf{q}^k,\mathbf{v})_D + (T^k,\nabla\cdot\mathbf{v})_D - \langle T_D,\mathbf{v}\cdot\mathbf{n}\rangle_{\partial D}, \label{eq:res1} \\
R_2(\varphi) &= -(\nabla\cdot\mathbf{q}^k,\varphi)_D + \beta\sigma\left(\sum_{m=1}^{N_a} w_m \Phi_m^k - (T^k)^4, \varphi\right)_D, \label{eq:res2} \\
R_{3,m}(\psi) &= (f_m,\psi)_D - (\Omega_m\cdot\nabla\Phi_m^k,\psi)_D - \sigma(\Phi_m^k,\psi)_D + \sigma((T^k)^4,\psi)_D. \label{eq:res3}
\end{align}
The boundary conditions for $\delta\Phi_m$ on inflow boundaries and for $\delta T$ on $\partial D$ are homogeneous, imposed as in \eqref{eq:boundary}. Here $\Phi^k = \sum_{m=1}^{N_a} w_m \Phi_m^k$ is the angular average of the intensity at the current iterate.

The well-posedness of the mixed formulation \eqref{eq:cont1}-\eqref{eq:cont3} will be investigated in a separate work,
which will not be addressed in the present work.
This continuous formulation serves as the point of departure for the spatial discretisation presented next.

\subsection{Mixed and DG upwind discretization for coupled conduction-radiation equation}
Let $\mathcal{T}_h$ be a quasi-uniform and shape-regular triangulation of $D$, and denote its elements by $K$. Let $\mathbb{P}_p(K)$ be the space of polynomials of degree at most $p$ on $K$. The corresponding vector-valued space is denoted by $\mathbb{P}_p(K)^3$.

The heat flux \(\mathbf{q}\) is discretised in the Raviart-Thomas finite element space
\[
\mathbf{W}_h = \bigl\{ \mathbf{v}_h \in \mathbf{H}(\mathrm{div},D) : \mathbf{v}_h|_K \in \mathbb{P}_p(K)^3 + \mathbf{x}\,\widetilde{\mathbb{P}}_p(K),\; \forall K\in\mathcal{T}_h \bigr\},
\]
where $\widetilde{\mathbb{P}}_p(K)$ denotes the space of homogeneous polynomials of degree \(p\).
The temperature $T$ is approximated in the discontinuous finite element space
\[
V_h = \bigl\{ v_h \in L^2(D) : v_h|_K \in \mathbb{P}_p(K),\; \forall K\in\mathcal{T}_h \bigr\},
\]
with no continuity enforced across element interfaces. We take the same polynomial degree $p$ for the temperature space to satisfy the discrete inf-sup condition for the mixed formulation~\cite{Brezzi1974, BrezziFortin1991, RaviartThomas1977}.

The radiative intensity $\mathbf{\Phi} = (\Phi_{1},\dots,\Phi_{N_a})^\top$  is approximated in the vector-valued discontinuous finite element space
\[
\mathbf{V}_h := \{\boldsymbol{\psi}_h\in[L^2(D)]^{N_a} : \boldsymbol{\psi}_h|_K \in [\mathbb{P}_p(K)]^{N_a},\; \forall K\in\mathcal{T}_h\},
\]
The DG discretisation of \eqref{eq:sd3} requires a numerical flux for the convective term.
On an interior face $e = \partial K^+\cap\partial K^-$, the upwind flux is defined by
\[
\widehat{\delta\Phi}_{m,h}^{\text{upw}} =
\begin{cases}
\delta\Phi_{m,h}|_{K^+}, & \Omega_m\cdot\mathbf{n}^+ > 0,\\
\delta\Phi_{m,h}|_{K^-}, & \Omega_m\cdot\mathbf{n}^+ < 0,
\end{cases}
\]
where $\mathbf{n}^+$ is the outward unit normal to $K^+$.
On the inflow boundary $\Gamma_{\mathrm{in}}^m = \{\mathbf{x}\in\partial D:\Omega_m\cdot\mathbf{n}<0\}$, the boundary condition gives $\widehat{\delta\Phi}_{m,h}^{\text{upw}} = 0$ due to the homogeneous increments.

For simplicity, we denote $\boldsymbol{\Phi}_h^k = (\Phi_{1,h}^k,\dots,\Phi_{N_a,h}^k)^\top$ and $\delta\boldsymbol{\Phi}_h = (\delta\Phi_{1,h},\dots,\delta\Phi_{N_a,h})^\top$. Given the current Newton iterate \((\mathbf{q}_h^k, T_h^k, \boldsymbol{\Phi}_h^k) \in \mathbf{W}_h\times V_h\times\mathbf{V}_h\), the fully discrete variational formulation reads: find $(\delta\mathbf{q}_h,\delta T_h,\delta\boldsymbol{\Phi}_h) \in \mathbf{W}_h\times V_h\times\mathbf{V}_h$ such that
\begin{align}
  (\kappa^{-1}\delta\mathbf{q}_h,\mathbf{v}_h)_D - (\delta T_h,\nabla\!\cdot\!\mathbf{v}_h)_D
    &= \mathcal{R}_1^{(k)}(\mathbf{v}_h) && \forall\mathbf{v}_h\in\mathbf{W}_h, \label{eq:dginc1}\\
  (\nabla\cdot\delta\mathbf{q}_h,\varphi_h)_D
    + 4\beta\sigma\bigl((T_h^k)^3\delta T_h,\varphi_h\bigr)_D
    - \beta\sigma\sum_{m=1}^{N_a} w_m (\delta\Phi_{m,h},\varphi_h)_D
    &= \mathcal{R}_2^{(k)}(\varphi_h) && \forall\varphi_h\in V_h, \label{eq:dginc2}\\
    \mathcal{F}(\Omega_m, \delta\Phi_{m,h}, \psi_h) - 4\sigma \sum_{K} \bigl((T_h^k)^3\delta T_h,\psi_h\bigr)_K + \sigma \sum_{K} \bigl(\delta\Phi_{m,h},\psi_h\bigr)_K &= \mathcal{R}_{3,m}^{(k)}(\psi_h) && \forall\psi_h\in V_h,
      \label{eq:dginc3}
\end{align}
where $m=1,\dots,N_a$, the upwind DG discretization of the advection term and residual functionals are defined by
\begin{align}
  \mathcal{F}(\Omega_m, \delta\Phi_{m,h}, \psi_h) &= -\sum_{K} \bigl( \delta\Phi_{m,h},\,\Omega_m\!\cdot\!\nabla\psi_h \bigr)_K
    + \sum_{K} \langle \widehat{\delta\Phi}_{m,h}^{\text{upw}},\,\psi_h \rangle_{\partial K,\Omega_m},\label{advection}\\
  \mathcal{R}_1^{(k)}(\mathbf{v}_h) &= -(\kappa^{-1}\mathbf{q}_h^k,\mathbf{v}_h)_D
    + (T_h^k,\nabla\!\cdot\!\mathbf{v}_h)_D
    - \langle T_D,\mathbf{v}_h\!\cdot\!\mathbf{n}\rangle_{\partial D}, \label{eq:calR1}\\
  \mathcal{R}_2^{(k)}(\varphi_h) &=
    - (\nabla\!\cdot\!\mathbf{q}_h^k,\varphi_h)_D
    + \beta\sigma\Bigl(\sum_{m=1}^{N_a} w_m \Phi_{m,h}^k - (T_h^k)^4,\varphi_h\Bigr)_D, \label{eq:calR2}\\
  \mathcal{R}_{3,m}^{(k)}(\psi_h) &= \sum_{K} (f_m,\psi_h)_K
    +\sum_{K} \bigl( \Phi_{m,h}^k,\,\Omega_m\!\cdot\!\nabla\psi_h \bigr)_K
    - \sum_{K} \langle \widehat{\Phi}_{m,h}^{k,\text{upw}},\,\psi_h \rangle_{\partial K,\Omega_m}\nonumber\\
      &\quad- \sigma\sum_{K} \bigl(\Phi_{m,h}^k,\psi_h\bigr)_K
      + \sigma\sum_{K} \bigl((T_h^k)^4,\psi_h\bigr)_K. \label{eq:calR3}
\end{align}

In the implementation, the $N_a$ angular intensities are stored as a single vector-valued field in \(\mathbf{V}_h\), with the spatial degree of freedom as the outer index and the angular direction as the inner index, following the convention adopted in~\cite{Badri2018, Badri2019, Badri2018Vectorial,Li2021}.
The well-posedness of the discrete system~\eqref{eq:dginc1}-\eqref{eq:dginc3} will also be analyzed in separate future work.
In the present work, we mainly propose the fully discrete scheme and devise robust block preconditioning for the resulting
linear algebraic systems.

\section{Block preconditioning method}

Each Newton step requires the solution of the linearised large-scale discrete system \eqref{eq:dginc1}-\eqref{eq:dginc3}.
The coefficient matrix is nonsymmetric and indefinite, which can cause iterative solvers to converge slowly or even fail without efficient preconditioning~\cite{Badri2019, Warsa2004}. Hence, the development of efficient and robust preconditioners is essential to speed up the iterative solution. Here, the robustness means that
\begin{itemize}
  \item the iteration count is robust with respect to the mesh size or the number of degrees of freedom.
  \item the number of iterations is robust with respect to the the physical parameters,
  which can vary over a relative wide range.
\end{itemize}

\subsection{Block structure of the linearised system}

After assembling the global coefficient matrix, the discrete linear system \eqref{eq:dginc1}-\eqref{eq:dginc3} takes the form
\begin{equation}
  \mathcal{A}\,\delta\mathbf{x} = \mathbf{b},
  \label{eq:sys}
\end{equation}
where
\begin{align*}
    \mathcal{A} =
  \begin{bmatrix}
    M_\mathbf{q} & -D^\top & 0 \\
    D & M_T & S \\
    0 & C & B_\mathbf{\Phi}
  \end{bmatrix},\quad
  \delta\mathbf{x} =
  \begin{pmatrix}
    \delta\mathbf{q}_h \\
    \delta T_h \\
    \delta\boldsymbol{\Phi}_h
  \end{pmatrix},\quad
  \mathbf{b} =
  \begin{pmatrix}
    \mathbf{b}_\mathbf{q} \\
    \mathbf{b}_T \\
    \mathbf{b}_\mathbf{\Phi}
  \end{pmatrix}.
\end{align*}

Following the standard multi-component ordering used in vector-valued finite elements~\cite{Badri2018Vectorial},
we organise $\delta\boldsymbol{\Phi}_h$ with the mesh element as the outer index and the angular direction as the inner index,
so that each mesh element carries $N_\mathbf{\Phi}$ basis functions, each with $N_a$ angular components, where $N_\mathbf{\Phi} = \dim \mathbb{P}_\ell(K)$ in the DG discretisation.
Applying this ordering, $B_\mathbf{\Phi}$ inherits the sparse structure
of the DG operator with $(N_\mathbf{\Phi} N_a)\times(N_\mathbf{\Phi} N_a)$ blocks.
As a result, the global matrix can be assembled with just one loop over the mesh.
At the same time, the coupling blocks $S$ and $C$ are assembled directly from the angular quadrature weights,
making the implementation both simple and efficient.

The blocks are summarised in Table~\ref{tab:blocks}.  The matrix $\mathcal{A}$ is nearly block upper-triangular: the only deviations from a strict block upper-triangular form are the coupling blocks $D$ and $C$ in the $(2,1)$ and $(3,2)$ positions. This structure motivates the construction of preconditioners based on block Gaussian elimination.

Let $\{\mathbf{v}_i : 1\leq i \leq N_q\}$ and $\{\varphi_i : 1\leq i \leq N_T\}$ be the bases of $\mathbf{W}_h$ and $V_h$, respectively.
For the radiation space $\mathbf{V}_h$, due to our vector finite element index order, the bases of $\mathbf{V}_h$ can be written as
$\{\boldsymbol{\psi}_{i,m} : 1\leq i\leq N_\mathbf{\Phi},\; 1\leq m\leq N_a\},$
where $\boldsymbol{\psi}_{i,m}$ has the $m$-th angular component equal to $\varphi_i$ and all other components zero.
Then the entries of the block matrices in \eqref{eq:sys} are defined as follows:
\begin{align*}
  (M_\mathbf{q})_{ij} &= ( \kappa^{-1} \mathbf{v}_j, \mathbf{v}_i )_D,
    & 1\leq i,j\leq N_q, \\
  (D)_{ij} &= ( \nabla\cdot \mathbf{v}_j, \varphi_i )_D,
    & 1\leq i\leq N_T,\;1\leq j\leq N_q, \\
  (M_T)_{ij} &= 4\beta\sigma ( (T_h^k)^3 \varphi_j, \varphi_i )_D,
    & 1\leq i,j\leq N_T, \\
  (S)_{i,(j,m)} &= -\beta\sigma w_m ( \varphi_j, \varphi_i )_D,
    & 1\leq i, j\leq N_T,\;1\leq m\leq N_a, \\
  (C)_{(i,m),j} &= -4\sigma ( (T_h^k)^3 \varphi_j, \varphi_i )_D,
    & 1\leq i, j\leq N_T,\;1\leq m\leq N_a, \\
  (B_\Phi)_{(i,m),(j,n)} &= \mathcal{F}(\Omega_m, \varphi_j, \varphi_i) + \sigma \sum_{K} (\varphi_j, \varphi_i)_K - \beta\sigma w_n \sum_{K} (\varphi_j, \varphi_i)_K,& 1\leq i,j\leq N_T,\;1\leq m,n\leq N_a,
\end{align*}
where $(i,m)$ corresponds to the linear index $i \cdot N_a + m$ (with indices counted from 0),
and $\mathcal{F}(\Omega_m, \cdot, \cdot)$ is the DG bilinear form defined in \eqref{advection}.
which generally couples the angular directions $m$ and $n$ through the upwind flux and the absorption term.

\begin{table}[h]
\centering
\caption{Block structure of the linearised system.}
\label{tab:blocks}
\begin{tabular*}{\textwidth}{@{} lll @{}}
\toprule
Block & Dimension & Role and structure \\
\midrule
$M_\mathbf{q}$ & $N_q\times N_q$ & weighted mass matrix; symmetric positive definite; sparse \\
$D$ & $N_T\times N_q$ & discrete divergence operator (from conduction equation); sparse \\
$M_T$ & $N_T\times N_T$ & reactive temperature mass matrix; diagonal for piecewise constant discretisation \\
$S$ & $N_T\times (N_\mathbf{\Phi} N_a)$ & temperature-to-radiation coupling; wide rectangular; diagonal in space \\
$C$ & $(N_\mathbf{\Phi} N_a)\times N_T$ & radiation-to-temperature coupling; tall rectangular; diagonal in space \\
$B_\Phi$ & $(N_\mathbf{\Phi} N_a)\times (N_\mathbf{\Phi} N_a)$ & upwind DG transport operator with angular coupling \\
\bottomrule
\end{tabular*}
\end{table}

\subsection{Block Schur complement preconditioner}

Block preconditioning method based on approximate Schur complements is a well-established technique for coupled multiphysics problems, where the coupling between different physical fields is the primary source of iterative stiffness. The general theory of preconditioning for saddle-point systems~\cite{Benzi2005} provides a framework for constructing such preconditioners. For the incompressible Navier-Stokes equations, Elman et al.~\cite{ElmanSilvesterWathen2014} demonstrated that Schur-complement-based preconditioners can achieve mesh-independent convergence by properly approximating the pressure Schur complement. For the stationary incompressible MHD equations, Li and Zheng~\cite{LiZheng2017} extended this idea to problems where the velocity-pressure and magnetic field blocks are coupled through the Lorentz force. Their strategies proves robust with respect to the mesh size and relevant physical parameters.

In the present coupled heat radiation transport system, the matrix $\mathcal{A}$ exhibits a similar block structure. We therefore construct a block upper-triangular preconditioner by eliminating the heat flux and temperature blocks
through approximate Schur complements. The approximate Schur complement for the temperature block, obtained by eliminating the heat flux, takes the form
\begin{equation}
  L_T := M_T + D\,\operatorname{diag}(M_\mathbf{q})^{-1}D^\top,
  \label{eq:LTdef}
\end{equation}
where the dense inverse of $M_\mathbf{q}$ is replaced by its diagonal approximation. This is standard practice in $H(\mathrm{div})$-based mixed methods, as the full inverse would be prohibitively expensive to form and apply, and is justified by the diagonal dominance of the mass matrix. For $P_{\mathrm{Schur}}$, we further eliminate the temperature-radiation coupling to obtain the radiation Schur complement
\begin{equation}
  L_\mathbf{\Phi} := B_\mathbf{\Phi} - C\,\operatorname{diag}(L_T)^{-1}S.
  \label{eq:Lphidef}
\end{equation}
The resulting block upper-triangular preconditioner is
\begin{equation}
  P_{\mathrm{Schur}} :=
  \begin{bmatrix}
    M_\mathbf{q} & -D^\top & 0 \\[4pt]
    0 & L_T & S \\[4pt]
    0 & 0 & L_\mathbf{\Phi}
  \end{bmatrix}.
\end{equation}
This preconditioner captures the dominant couplings between the heat flux, temperature, and radiation intensity, while its block triangular form allows for efficient application by back-substitution. The diagonal approximations avoid the expense of forming dense Schur complements, making the preconditioner practical for large-scale problems.

Since $P_{\mathrm{Schur}}$ is block upper-triangular, its application to a residual vector $(\mathbf{r}_\mathbf{q}, \mathbf{r}_T, \mathbf{r}_\mathbf{\Phi})^\top$ proceeds by block back-substitution:
\begin{enumerate}
  \item Solve $L_\mathbf{\Phi}\,\mathbf{x}_\mathbf{\Phi} = \mathbf{r}_\mathbf{\Phi}$.
  \item Solve $L_T\,\mathbf{x}_T = \mathbf{r}_T - S\,\mathbf{x}_\mathbf{\Phi}$.
  \item Solve $M_\mathbf{q}\,\mathbf{x}_\mathbf{q} = \mathbf{r}_\mathbf{q} + D^\top \mathbf{x}_T$.
\end{enumerate}
The matrices $L_T$ and $L_\mathbf{\Phi}$ are sparse and can be solved by suitable iterative solvers (e.g., AMG-preconditioned CG for $L_T$ and GMRES for $L_\mathbf{\Phi}$).

\subsection{Split preconditioner}

Operator splitting is a classical strategy for problems where different physical processes evolve on significantly different time scales~\cite{Mousseau2000, BrownWoodward2001}. In the present radiation-conduction system, the characteristic time scale of radiation transport is typically much shorter than that of thermal conduction. This motivates a split treatment of the two physics at the preconditioner level, analogous to explicit operator splitting at the time-stepping level, but applied algebraically within the preconditioner.

Recall that $L_\mathbf{\Phi} = B_\mathbf{\Phi} - C\,\operatorname{diag}(L_T)^{-1}S$ is the radiation Schur complement in $P_{\mathrm{Schur}}$, with $L_T = M_T + D\,\operatorname{diag}(M_\mathbf{q})^{-1}D^\top$. The term
\[
D\,\operatorname{diag}(M_{\mathbf{q}})^{-1}D^\top \;\Leftrightarrow\; -\nabla\cdot(\kappa\nabla), \quad \text{on } V_h
\]
is the discrete analogue of the diffusion operator. It arises naturally from the mixed formulation of the conduction equation after eliminating the heat flux via static condensation and represents the coupling between the temperature and heat flux variables.

Inspired by the explicit treatment of the diffusion term in transient operator-split schemes, we apply a similar idea at the algebraic level: we drop this diffusion coupling from the block $L_\mathbf{\Phi}$. This effectively decouples the temperature and radiation blocks. We therefore refer to this variant as the split preconditioner, denoted by $P_{\mathrm{Split}}$.

{When the diffusion and mass contributions are comparable, the effectiveness of $P_{\mathrm{Split}}$ is explained by the low-frequency structure of $L_T$. In Krylov methods, convergence is typically limited by the smallest eigenvalues of the preconditioned operator, which correspond to smooth error components that are difficult to eliminate. For such error components $v \in V_h$, the diffusion term satisfies
\[
D\,\operatorname{diag}(M_\mathbf{q})^{-1}D^\top v \approx 0,
\]
since constant vectors lie in its zero space and $v$ varies slowly. On the low-frequency subspace, this implies
\[
L_T v = (M_T + D\,\operatorname{diag}(M_\mathbf{q})^{-1}D^\top)v \approx M_T v,
\]
i.e. $L_T^{-1} v \approx M_T^{-1} v$. The split preconditioner thus provides an accurate approximation precisely on the modes that are hardest for Krylov methods to eliminate.}

When the coupled system is radiation-dominated ($\kappa\ll 1$), the diffusion coupling $D\,\operatorname{diag}(M_{\mathbf{q}})^{-1}D^\top$ is much weaker than the reactive term $M_T$. In this regime, replacing $L_T$ with $M_T$ is a reasonable approximation in $L_\mathbf{\Phi}$. Let $L_S:= B_\mathbf{\Phi} - C\,\operatorname{diag}(M_T)^{-1}S$, then the resulting preconditioner is
\begin{equation}
P_{\mathrm{Split}} :=
\begin{bmatrix}
  M_{\mathbf{q}} & -D^\top & 0 \\[4pt]
  0 & L_T & S \\[4pt]
  0 & 0 & L_S
\end{bmatrix}.
\label{eq:PSplit}
\end{equation}
This approximation is expected to be effective when the diffusion term of temperature equation
is weak relative to the mass term, while in conduction-dominated cases the $P_{\mathrm{Schur}}$ should be
recommended.

\subsection{Block Jacobi preconditioner}

Block Jacobi preconditioning is one of the simplest and most
widely used strategies~\cite{BrownWoodward2001, ElmanSilvesterWathen2014}. Its core idea is to approximate the system matrix by its diagonal blocks, neglecting the off-diagonal couplings that are less important. This strategy is particularly attractive when the coupling between blocks is weak relative to the dominant physics within each block.

In the present work, this regime occurs when $\sigma$ is large.
In the optically thick limit, the transport operator $B_\mathbf{\Phi}$ is strongly diagonally dominant, so $B_\mathbf{\Phi}$ alone provides a reasonable approximation to the full radiation Schur complement $L_\mathbf{\Phi}$. The radiative intensity is then essentially determined by the local temperature via the Planck source term, and the coupling through $S$ and $C$ becomes less important. This motivates the following block Jacobi preconditioner
\begin{equation}
  P_{\mathrm{BJ}} :=
  \begin{bmatrix}
    M_\mathbf{q} & -D^\top & 0 \\[4pt]
    0 & L_T & S \\[4pt]
    0 & 0 & B_\mathbf{\Phi}
  \end{bmatrix},
\end{equation}
where $B_\mathbf{\Phi}$ is the original radiation block without Schur-complement modification.

$P_{\mathrm{BJ}}$ is the simplest of the three preconditioners, replacing the radiation Schur complement by the original operator $B_\mathbf{\Phi}$. This is justified in the conduction-dominated regime, where the coupling between the temperature and radiation blocks is weak. In radiation-dominated regimes, where this coupling becomes strong, the $P_{\mathrm{Schur}}$ should be recommended.

\section{Numerical experiments}\label{sec:numerical}

In this section, we present numerical experiments that verify the analyses of the preceding sections. All computations are performed with the parallel finite element library \textsc{PHG}~\cite{Zhang2009} on a laptop equipped with an Intel Core i5-11300H processor (4 cores, 8 threads, 3.10\,GHz) and 16\,GB of RAM, running Ubuntu Linux under Windows Subsystem for Linux (WSL 2), using 4 MPI processes.

The spatial discretisation employs the mixed Raviart-Thomas-discontinuous Galerkin pair $\mathrm{RT}_k$-$\mathrm{DG}_k$ for the flux-temperature block, and the discontinuous Galerkin space $\mathrm{DG}_k$ for the radiation intensity $\mathbf{\Phi}$. The polynomial degree $k$ is given in each example. The angular discretisation employs the $S_2$ level symmetric quadrature with $N_a = 8$ directions. Without specifications, the initial guess is taken as $\mathbf{q}^0 = \mathbf{0}$, $T^0 = 1.0$, $\mathbf{\Phi}^0 = \mathbf{0}$. The Newton iteration is terminated when the relative tolerance falls below $\epsilon = 10^{-5}$.

At each Newton step the linearised system is solved by right-preconditioned FGMRES with a restart of 100, with a relative tolerance of $\varepsilon = 10^{-6}$. The inner solves required by the block preconditioners are performed as follows. The $M_\mathbf{q}$ block is solved by PCG with the diagonal preconditioner from hypre. The $L_T$ block is solved by PCG with BoomerAMG~\cite{HensonYang2002} from hypre. The radiation-related blocks $B_\mathbf{\Phi}$, $L_\mathbf{\Phi}$, and $L_S$ are solved by FGMRES~\cite{Saad1993} with the additive Schwarz preconditioner ASM(1) and ILU(2) from PETSc~\cite{Balay2025}. All inner solves in this example use a relative tolerance of $\varepsilon_0 = 10^{-2}$. For $P_{\mathrm{ASM}}$, the global system is solved by FGMRES with the same ASM(1) preconditioner with ILU(2) from PETSc.

\subsection{Example 1: Convergence order of the numerical scheme}

In this example, we verify the spatial convergence rates of the mixed finite element discretisation. The numerical convergence assessment begins with the lowest-order pair $\mathrm{RT}_0-\mathrm{DG}_0-\mathrm{DG}_0$ described in Section~3.  The physical parameters are taken as $\kappa = \beta = \sigma = 1$, and the computational domain is the unit cube $\Omega = [0,1]^3$.  The discrete ordinates quadrature uses $S_2$ with $N_a = 8$ directions.

The exact solution is prescribed via the method of manufactured solutions:
\begin{align*}
  \mathbf{q}(\mathbf{x}) &= -\kappa\,\nabla T
    = -\tfrac12\kappa\pi\cos(\pi(x+y+z))\,(1,1,1)^\top,\\[2pt]
  T(\mathbf{x}) &= 1 + \tfrac12\sin(\pi(x+y+z)),\\[2pt]
  \Phi(\mathbf{x},\Omega) &= \sin(\pi(x+y+z)) + x^2 + y^2 + z^2 + xyz.
\end{align*}
The radiative intensity is taken to be isotropic (independent of $\Omega$) so that the angular quadrature does not introduce additional errors beyond those inherent in the spatial discretisation, thereby isolating the spatial convergence behaviour.  The source terms in the governing equations are then obtained by inserting the manufactured solution into the equations.

\begin{table}[h]
\centering
\caption{DOF counts for the high-order discretisation $\mathrm{RT}_0-\mathrm{DG}_0-\mathrm{DG}_0$.}
\label{tab:mesh}
\begin{tabular}{@{} llllll @{}}
\toprule
Mesh & $h$ & DOF($\mathbf{q}_h$) & DOF($T_h$) & DOF($\boldsymbol{\Phi}_h$) & Total DOF \\
\midrule
$\mathcal{T}_1$ & $3.750\times10^{-1}$ & 1{,}632  & 768    & 6{,}144    & 8{,}544 \\
$\mathcal{T}_2$ & $1.875\times10^{-1}$ & 12{,}672 & 6{,}144& 49{,}152   & 67{,}968 \\
$\mathcal{T}_3$ & $9.375\times10^{-2}$ & 99{,}840 & 49{,}152& 393{,}216 & 542{,}208 \\
$\mathcal{T}_4$ & $4.688\times10^{-2}$ & 792{,}576&393{,}216&3{,}145{,}728&4{,}331{,}520 \\
\bottomrule
\end{tabular}
\end{table}

\begin{table}[h]
\centering
\caption{Errors and convergence rates for $\mathrm{RT}_0-\mathrm{DG}_0-\mathrm{DG}_0$.}
\label{tab:convergence}
\begin{tabular}{@{} llllllll @{}}
\toprule
Mesh & $h$ & $\|\mathbf{q}-\mathbf{q}_h\|_{0}$ & Order & $\|T-T_h\|_{0}$ & Order & $\|\boldsymbol{\Phi}-\boldsymbol{\Phi}_h\|_{0}$ & Order \\
\midrule
$\mathcal{T}_1$ & $3.750\times10^{-1}$ & $5.107\times10^{-1}$ & ---      & $7.961\times10^{-2}$ & --- & $8.256\times10^{-1}$ & --- \\
$\mathcal{T}_2$ & $1.875\times10^{-1}$ & $2.462\times10^{-1}$ & 1.05 & $3.912\times10^{-2}$ & 1.03 & $4.546\times10^{-1}$ & 0.86 \\
$\mathcal{T}_3$ & $9.375\times10^{-2}$ & $1.211\times10^{-1}$ & 1.02 & $1.945\times10^{-2}$ & 1.01 & $2.453\times10^{-1}$ & 0.89 \\
$\mathcal{T}_4$ & $4.688\times10^{-2}$ & $6.010\times10^{-2}$ & 1.01 & $9.708\times10^{-3}$ & 1.00 & $1.286\times10^{-1}$ & 0.93 \\
\bottomrule
\end{tabular}
\end{table}

Table~\ref{tab:mesh} shows the information for the meshes and the computation scale of the block and the global matrix. From Table~\ref{tab:convergence}, we find that the convergence rates for $(\mathbf{q}_h, T_h, \mathbf{\Phi}_h)$ are given by
\begin{align}
  ||\mathbf{q} - \mathbf{q}_h||_{\mathbf{0}} \sim \mathcal{O}(h),\quad
  ||T - T_h||_{\mathbf{0}} \sim \mathcal{O}(h),\quad
  ||\mathbf{\Phi}-\mathbf{\Phi}_h||_{\mathbf{0}}\sim\mathcal{O}(h).
\end{align}

Remember that we are using the lowest-order Raviart-Thomas face finite element for discretizing $\mathbf{q}$, the piecewise constant finite element for discretizing $T$ and $\mathbf{\Phi}_h$. This means that the optimal convergence rates are obtained for all variables.

The same manufactured solution and angular quadrature ($S_2$ with $N_a=8$) are used to test the higher-order pair $\mathrm{RT}_1-\mathrm{DG}_1-\mathrm{DG}_1$. Table~\ref{tab:mesh_ho} lists the degrees of freedom on the three mesh levels used for this test.

\begin{table}[h]
\centering
\caption{DOF counts for the high-order discretisation $\mathrm{RT}_1-\mathrm{DG}_1-\mathrm{DG}_1$.}
\label{tab:mesh_ho}
\begin{tabular}{@{} llllll @{}}
\toprule
Mesh & $h$ & DOF($\mathbf{q}_h$) & DOF($T_h$) & DOF($\boldsymbol{\Phi}_h$) & Total DOF \\
\midrule
$\mathcal{G}_1$ & $2.795\times10^{-1}$ &   14{,}400 &   6{,}144 &   49{,}152 &   69{,}696 \\
$\mathcal{G}_2$ & $1.398\times10^{-1}$ &  112{,}896 &  49{,}152 &  393{,}216 &  555{,}264 \\
$\mathcal{G}_3$ & $6.988\times10^{-2}$ &  893{,}952 & 393{,}216 & 3{,}145{,}728 & 4{,}432{,}896 \\
\bottomrule
\end{tabular}
\end{table}

\begin{table}[h]
\centering
\caption{Errors and convergence rates for $\mathrm{RT}_1-\mathrm{DG}_1-\mathrm{DG}_1$.}
\label{tab:conv_ho}
\begin{tabular}{@{} llllllll @{}}
\toprule
Mesh & $h$ & $\|\mathbf{q}-\mathbf{q}_h\|_{0}$ & Order & $\|T-T_h\|_{0}$ & Order & $\|\boldsymbol{\Phi}-\boldsymbol{\Phi}_h\|_{0}$ & Order \\
\midrule
$\mathcal{G}_1$ & $2.795\times10^{-1}$ & $7.433\times10^{-2}$ & ---  & $1.261\times10^{-2}$ & --- & $1.188\times10^{-1}$ & --- \\
$\mathcal{G}_2$ & $1.398\times10^{-1}$ & $1.877\times10^{-2}$ & 1.99 & $3.212\times10^{-3}$ & 1.97 & $3.073\times10^{-2}$ & 1.95 \\
$\mathcal{G}_3$ & $6.988\times10^{-2}$ & $4.698\times10^{-3}$ & 2.00 & $8.066\times10^{-4}$ & 1.99 & $7.761\times10^{-3}$ & 1.98 \\
\bottomrule
\end{tabular}
\end{table}

Table~\ref{tab:conv_ho} shows the information for the meshes and the computation scale of the block and the global matrix. From Table~\ref{tab:conv_ho}, we find that the convergence rates for $(\mathbf{q}_h, T_h, \mathbf{\Phi}_h)$ are given by
\begin{align}
  ||\mathbf{q} - \mathbf{q}_h||_{\mathbf{0}} \sim \mathcal{O}(h^2),\quad
  ||T - T_h||_{\mathbf{0}} \sim \mathcal{O}(h^2),\quad
  ||\mathbf{\Phi}-\mathbf{\Phi}_h||_{\mathbf{0}}\sim\mathcal{O}(h^2).
\end{align}

Remember that we are using the first-order Raviart-Thomas face finite element for discretizing $\mathbf{q}$, the piecewise constant finite for discretizing $T$ and $\mathbf{\Phi}_h$. This means that the optimal convergence rates are obtained for all variables.

\subsection{Example 2: Performance of block preconditioners}

This subsection examines the robustness of the three block preconditioners from Section~4 through a comparison with the additive Schwarz method $(P_{\mathrm{ASM}}$) from PETSc. All tests use the lowest-order discretisation $\mathrm{RT}_0-\mathrm{DG}_0-\mathrm{DG}_0$ described in Section~3. The absorption coefficient is taken to be strongly discontinuous with four alternating layers along the $z$-direction:
\[
\sigma(z) =
\begin{cases}
100,   & 0.00 \le z < 0.25,\; 0.50 \le z < 0.75,\\[4pt]
0.01,  & 0.25 \le z < 0.50,\; 0.75 \le z \le 1.00.
\end{cases}
\]
The large contrast in $\sigma$ is intended to test the robustness of the preconditioners under strongly heterogeneous material properties. We consider three test configurations that differ by the thermal conductivity $\kappa$:

\begin{itemize}
  \item Case~1 (conduction-dominated): $\kappa = 10.0$, $\beta = 1$.  With large thermal conductivity, heat conduction dominates the energy transfer and the coupling to radiation is relatively weak.
  \item Case~2 (mixed): $\kappa = 10^{-2}$, $\beta = 1$. At this intermediate value, neither conduction nor radiation is negligible, and both mechanisms contribute to the energy balance.
  \item Case~3 (radiation-dominated): $\kappa = 10^{-5}$, $\beta = 1$.  With very small thermal conductivity, radiation dominates the energy transfer, leading to strong coupling between the temperature and radiation fields.
\end{itemize}

As will be shown in the following, all preconditioners maintain stable iteration counts despite the large jump in $\sigma$.
\begin{table}[h]
\centering
\small
\caption{FGMRES iteration counts for the three test configurations. Entries are given as $N_{\mathrm{FGMRES}}$ ($N_{\mathrm{Newton}}$).}
\label{tab:preconditioner_performance}
\begin{tabular}{l | c c c | c c c | c c c}
\toprule
\multirow{2}{*}{Preconditioner} & \multicolumn{3}{c}{Case 1: $\kappa=10$} & \multicolumn{3}{c}{Case 2: $\kappa=10^{-2}$} & \multicolumn{3}{c}{Case 3: $\kappa=10^{-5}$} \\
\cmidrule(lr){2-4} \cmidrule(lr){5-7} \cmidrule(lr){8-10}
 & $\mathcal{M}_1$ & $\mathcal{M}_2$ & $\mathcal{M}_3$ & $\mathcal{M}_1$ & $\mathcal{M}_2$ & $\mathcal{M}_3$ & $\mathcal{M}_1$ & $\mathcal{M}_2$ & $\mathcal{M}_3$ \\
\midrule
$P_{\mathrm{Schur}}$ & 11(4) & 11(3) & 10(3) & 20(6) & 35(6) & 49(5) & 6(6) & 6(6) & 9(6) \\
$P_{\mathrm{Split}}$ & 49(4) & 67(3) & 90(3) & 14(6) & 17(6) & 21(5) & 6(6) & 7(6) & 9(6) \\
$P_{\mathrm{BJ}}$ & 11(4) & 11(3) & 10(3) & 56(6) & 66(6) & 69(5) & 43(6) & 51(6) & 61(6) \\
$P_{\mathrm{ASM}}$ & 18(3) & 30(3) & 57(3) & 20(6) & 40(6) & 80(5) & 20(6) & 34(6) & 65(6) \\
\bottomrule
\end{tabular}
\end{table}

Table~\ref{tab:preconditioner_performance} reports the FGMRES iteration counts for the three test configurations across meshes $\mathcal{M}_1$, $\mathcal{M}_2$, and $\mathcal{M}_3$ ($3.4\times10^4$, $2.7\times10^5$, and $2.2\times10^6$ DOFs). Each entry is given as $N_{\mathrm{FGMRES}}$ ($N_{\mathrm{Newton}}$), where $N_{\mathrm{FGMRES}}$ denotes the average number of GMRES iterations and $N_{\mathrm{Newton}}$ is the number of Newton iterations.

The results reveal a clear regime-dependent hierarchy. In the conduction-dominated case ($\kappa=10$), $P_{\mathrm{Schur}}$ and $P_{\mathrm{BJ}}$ perform nearly identically and are mesh-independent, while $P_{\mathrm{ASM}}$ degrades with mesh refinement. In the mixed regime ($\kappa=10^{-2}$), $P_{\mathrm{Split}}$ substantially outperforms the others. In the radiation-dominated case ($\kappa=10^{-5}$), $P_{\mathrm{Schur}}$ and $P_{\mathrm{Split}}$ both perform well, whereas $P_{\mathrm{BJ}}$ deteriorates markedly. Overall, the iteration counts remain stable under mesh refinement, confirming the mesh robustness of the proposed preconditioners.

\begin{figure}[h]
	\centering
	\includegraphics[width=1\textwidth]{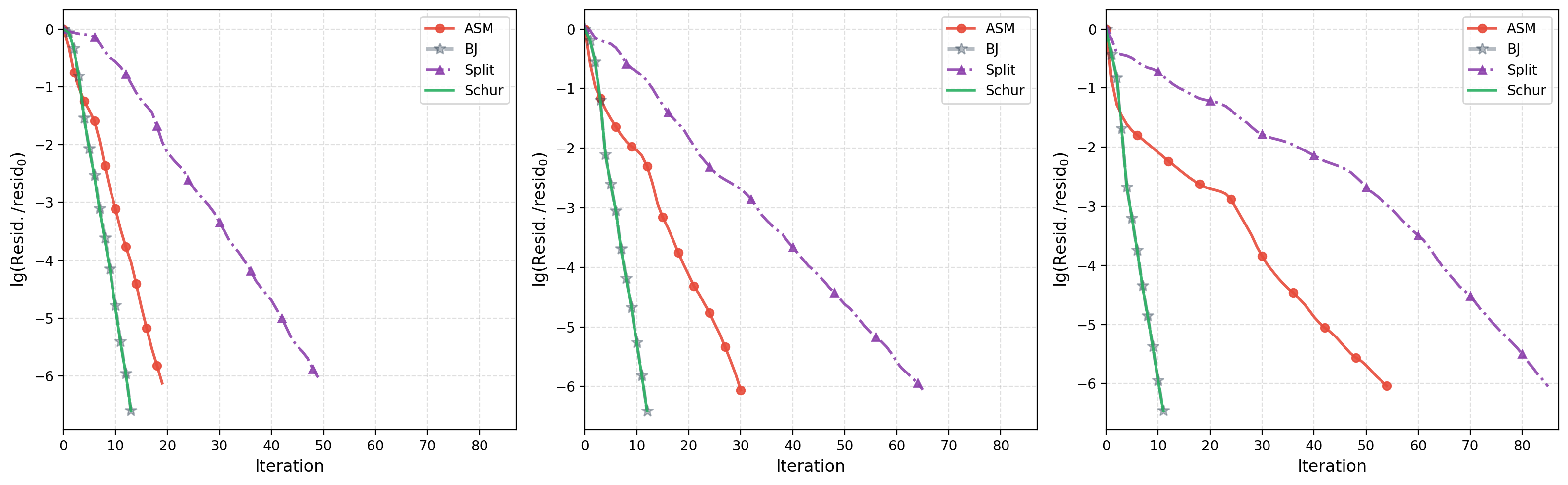}
	\caption{FGMRES convergence history for the four preconditioners on mesh $\mathcal{M}_1\sim\mathcal{M}_3$ (Case~1, $\kappa = 10.0$).}
	\label{fig:case1}
\end{figure}

The convergence histories are consistent with the iteration counts reported in Table~\ref{tab:preconditioner_performance}. Figure~\ref{fig:case1} shows the FGMRES convergence at the first Newton step on mesh $\mathcal{M}_1\sim\mathcal{M}_3$ for Case~1. $P_{\mathrm{Schur}}$ and $P_{\mathrm{BJ}}$ converge rapidly, whereas $P_{\mathrm{split}}$ and $P_{\mathrm{ASM}}$ require substantially more iterations.

\begin{figure}[h]
	\centering
	\includegraphics[width=1\textwidth]{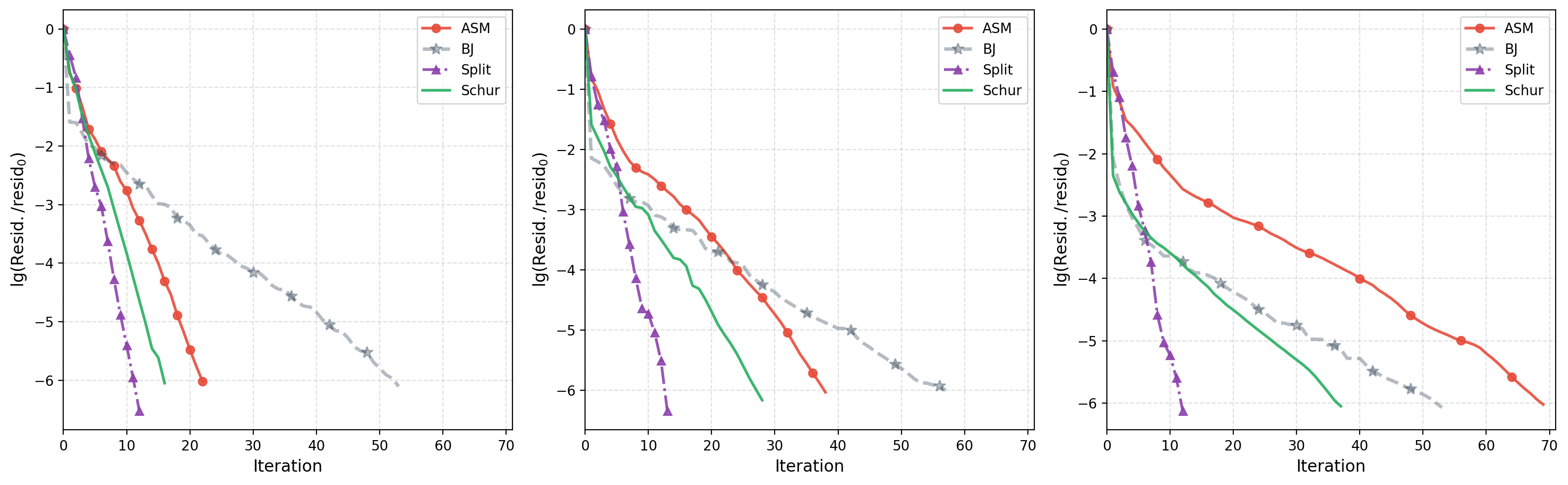}
	\caption{FGMRES convergence history for the four preconditioners on mesh $\mathcal{M}_1\sim\mathcal{M}_3$ (Case~2, $\kappa = 1e-2$).}
	\label{fig:case2}
\end{figure}

{The degradation of $P_{\mathrm{Schur}}$ in the mixed regime is explained by the strong coupling between blocks
\[
S \sim \sigma, \qquad C \sim \sigma.
\]
In the mixed regime, the diagonal approximation of $L_T$ is no longer accurate. The resulting error in the radiation Schur complement satisfies
\[
\|L_\Phi - \tilde{L}_\Phi\| = \|C\,\bigl(L_T^{-1} - \operatorname{diag}(L_T)^{-1}\bigr)\,S\| \sim \sigma^2\|L_T^{-1} - \operatorname{diag}(L_T)^{-1}\|.
\]
Consequently, the coupling blocks compound the diagonal approximation error as $\sigma$ increases, which explains the degraded performance of $P_{\mathrm{Schur}}$ in the mixed regime. The split preconditioner avoids this issue by dropping the diagonal approximation of $L_T$ entirely. In the conduction- and radiation-dominated regimes the diagonal approximation of $L_T$ is accurate, so this amplification does not cause performance degradation. In all cases, iteration counts remain stable under mesh refinement, confirming the mesh robustness of all proposed preconditioners.}

\begin{figure}[h]
	\centering
	\includegraphics[width=1\textwidth]{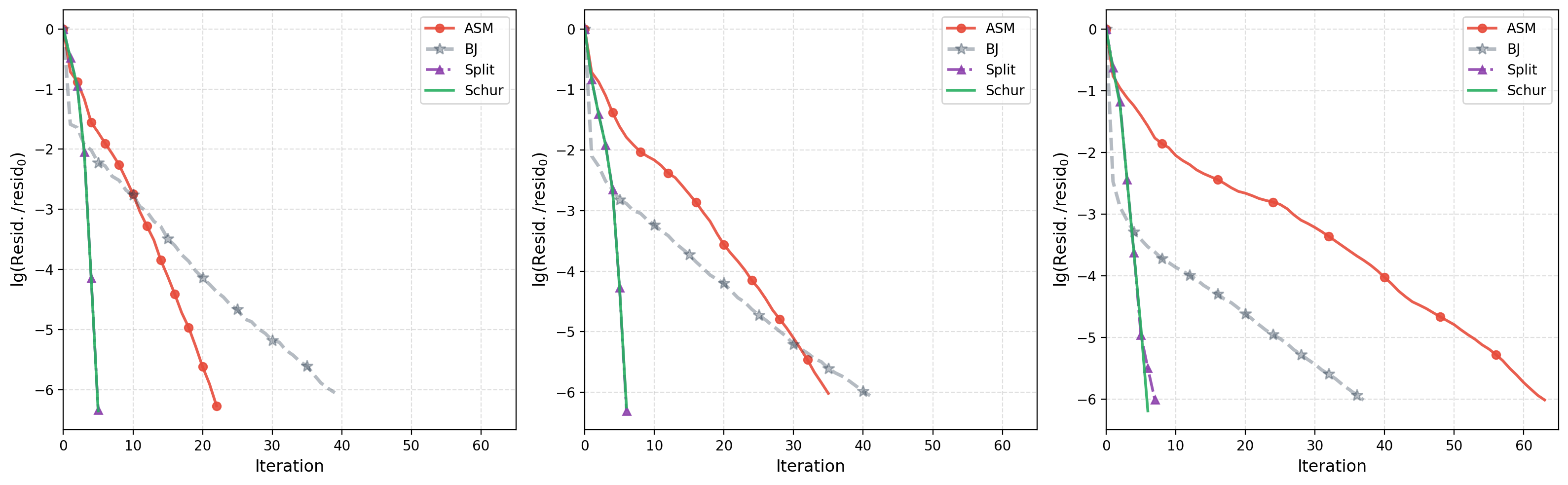}
	\caption{FGMRES convergence history for the four preconditioners on mesh $\mathcal{M}_1\sim\mathcal{M}_3$ (Case~3, $\kappa = 1e-5$).}
	\label{fig:case3}
\end{figure}

The convergence histories for Case~2 and Case~3 are shown in Figures~\ref{fig:case2} and~\ref{fig:case3}, respectively. In the mixed regime (Case~2), $P_{\mathrm{Split}}$ converges most rapidly, with $P_{\mathrm{Schur}}$ also performing well. In the radiation-dominated regime (Case~3), both $P_{\mathrm{Schur}}$ and $P_{\mathrm{Split}}$ converge rapidly. In both cases, $P_{\mathrm{BJ}}$ and $P_{\mathrm{ASM}}$ are significantly slower.

To summarise the hierarchy of the three preconditioners, Table~\ref{tab:hierarchy_summary} collects the iteration counts on mesh $\mathcal{M}_3$ for all configurations.
\begin{table}[h]
\centering
\small
\caption{Summary of the hierarchy of approximations: iteration counts on mesh $\mathcal{M}_3$ across the three test configurations.}
\label{tab:hierarchy_summary}
\begin{tabular}{l c c c}
\toprule
Preconditioner & Case 1: $\kappa=10$ & Case 2: $\kappa=10^{-2}$ & Case 3: $\kappa=10^{-5}$ \\
\midrule
$P_{\mathrm{Schur}}$ (full coupling) & 10(3) & 49(5) & 9(6) \\
$P_{\mathrm{Split}}$ (diffusion omitted) & 90(3) & 21(5) & 9(6) \\
$P_{\mathrm{BJ}}$ (coupling neglected) & 10(3) & 69(5) & 61(6) \\
\bottomrule
\end{tabular}
\end{table}

{The hierarchy of approximations is clearly reflected in the results. $P_{\mathrm{Schur}}$ is the most robust preconditioner overall, performing excellently in the conduction- and radiation-dominated regimes. In the mixed regime, its performance remains acceptable, with only a modest impact from the diagonal approximation error of $L_T$. $P_{\mathrm{Split}}$ is most efficient in the mixed regime, and performs comparably to $P_{\mathrm{Schur}}$ in the radiation-dominated case. $P_{\mathrm{BJ}}$ is competitive only in the conduction-dominated regime, where it matches the performance of $P_{\mathrm{Schur}}$.}

The sensitivity of the results to the inner-solver tolerance $\varepsilon_0$ is examined in the conduction-dominated regime ($\kappa=10$, $\beta=1$). With mesh $\mathcal{T}_4$ fixed, $\varepsilon_0$ is varied from $10^{-2}$ to $10^{-4}$. Table~\ref{tab:inner_tol} reports the average FGMRES iteration counts for the four preconditioners.

\begin{table}[h]
\centering
\small
\caption{Iteration counts for different inner solver tolerances in the conduction-dominated regime ($\kappa=10$, mesh $\mathcal{T}_4$). Entries are given as $N_{\mathrm{GMRES}}$ ($N_{\mathrm{Newton}}$).}
\label{tab:inner_tol}
\begin{tabular}{l c c c}
\toprule
Preconditioner & $\varepsilon_0=10^{-2}$ & $\varepsilon_0=10^{-3}$ & $\varepsilon_0=10^{-4}$ \\
\midrule
$P_{\mathrm{Schur}}$  & 10(3) & 10(3) & 10(3) \\
$P_{\mathrm{Split}}$  & 91(3) & 63(3) & 59(3) \\
$P_{\mathrm{BJ}}$     & 10(3) & 10(3) & 10(3) \\
$P_{\mathrm{ASM}}$    & 80(5) & 80(5) & 80(5) \\
\bottomrule
\end{tabular}
\end{table}

The tolerance $\varepsilon_0$ has little influence on the performance of $P_{\mathrm{Schur}}$ and $P_{\mathrm{BJ}}$, both of which maintain constant iteration counts. $P_{\mathrm{ASM}}$ also shows no sensitivity, as it has no sub-block solves. Only $P_{\mathrm{Split}}$ is affected: its iteration count decreases as the tolerance is tightened from $10^{-2}$ to $10^{-4}$, with the improvement beyond $10^{-2}$ being marginal. Based on these observations, $\varepsilon_0 = 10^{-2}$ offers a reasonable balance between accuracy and efficiency, and is used as the default relative tolerance for all inner solves in this paper.

\subsection{Example 3: A cubic cavity benchmark}\label{sec:example3}

To assess the performance of the proposed preconditioner on a physically realistic problem, we consider a three-dimensional cubic enclosure with a gray absorbing-emitting medium, following the benchmark configuration of Talukdar et al.~\cite{Talukdar2008}. The domain is the unit cube $\Omega = [0,1]^3$. The bottom wall ($z = 0$) is held at $T = 1$, while the other five walls are at $T = 0.5$. All walls are black, so the inflow intensity is given by $\Phi^{\mathrm{in}} = T_{\mathrm{wall}}^4$.

\begin{table}[h]
\centering
\renewcommand{\arraystretch}{1.4}
\caption{Average FGMRES iterations per Newton step
           $N_{\mathrm{GMRES}}(N_{\mathrm{Newton}})$ for $P_{\mathrm{Schur}}$(left) and $P_{\mathrm{ASM}}$(right).}
\begin{minipage}[t]{0.46\textwidth}
  \centering
  \label{tab:cavity_schur}
  \begin{tabular*}{\linewidth}{@{\extracolsep{\fill}} llll @{}}
  \toprule
  $\kappa$ & $\beta = 1.0$ & $\beta = 10.0$ & $\beta = 100.0$ \\
  \midrule
  $0.01$  &  $12.3(3)$  &  $11.4(5)$  &  $11.0(10)$ \\
  $1.0$   &  $9.7(3)$   &  $9.7(3)$   &  $9.8(4)$ \\
  $100.0$ &  $11.3(3)$  &  $11.3(3)$  &  $11.3(3)$ \\
  \bottomrule
  \end{tabular*}
\end{minipage}
\hfill
\begin{minipage}[t]{0.46\textwidth}
  \centering
  \label{tab:cavity_asm}
  \begin{tabular*}{\linewidth}{@{\extracolsep{\fill}} llll @{}}
  \toprule
  $\kappa$ & $\beta = 1.0$ & $\beta = 10.0$ & $\beta = 100.0$ \\
  \midrule
  $0.01$  &  $76.3(3)$  &  $53.6(5)$  &  $40.8(10)$ \\
  $1.0$   &  $79.0(3)$  &  $82.0(3)$  &  $67.8(4)$ \\
  $100.0$ &  $75.0(3)$  &  $77.0(3)$  &  $78.7(3)$ \\
  \bottomrule
  \end{tabular*}
\end{minipage}
\end{table}

The external source is set to zero ($f=0$) and the absorption coefficient is taken as $\sigma = 1.0$. The simulations are performed using the lowest-order discretisation $\mathrm{RT}_0-\mathrm{DG}_0-\mathrm{DG}_0$ on mesh $\mathcal{T}_4$. Table~\ref{tab:cavity_schur} reports the average FGMRES iteration counts for $P_{\mathrm{Schur}}$ and $P_{\mathrm{ASM}}$ across nine combinations of $\kappa$ and $\beta$. The computed temperature profiles are found to be in good agreement with the benchmark results of Talukdar et al.~\cite{Talukdar2008}.

Two observations from Table~\ref{tab:cavity_schur} are noteworthy. $P_{\mathrm{Schur}}$ outperforms $P_{\mathrm{ASM}}$ by approximately an order of magnitude across all parameter combinations, while its iteration count remains essentially constant as $\kappa$ and $\beta$ vary over a wide range. In contrast, $P_{\mathrm{ASM}}$ shows clear parameter dependence. This confirms that the block-structured preconditioner is robust on physically realistic configurations with non-trivial boundary conditions.

\begin{figure}[h]
	\centering
	\includegraphics[width=.45\textwidth]{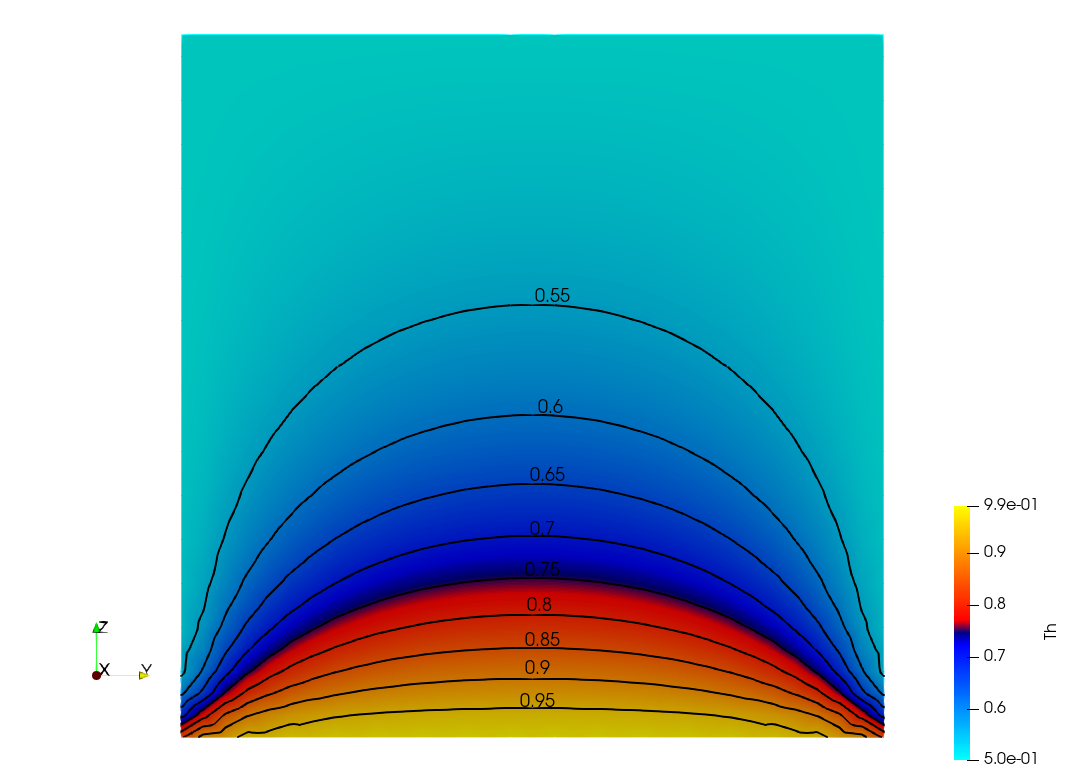}\hfill
	\includegraphics[width=.45\textwidth]{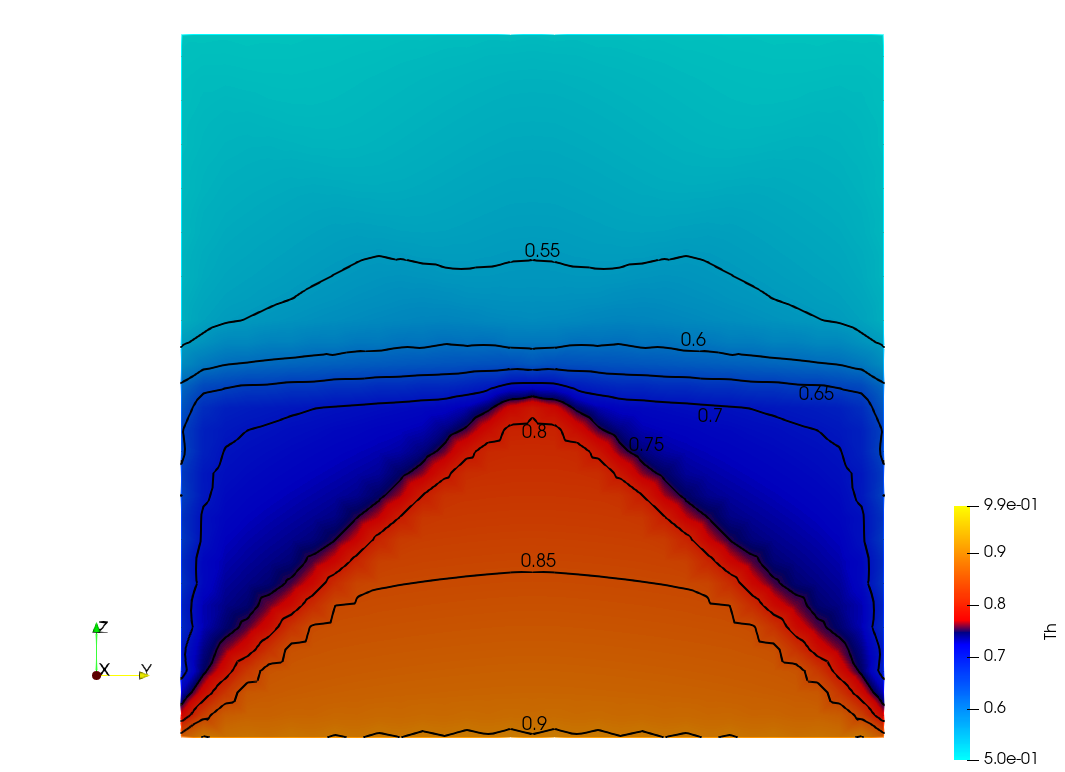}
	\caption{Temperature field and isotherms on $x = 0.5$.  Left: conduction-dominated ($\kappa = 100.0$, $\beta = 1.0$).  Right: radiation-dominated ($\kappa = 0.01$, $\beta = 100.0$).}
	\label{fig:cavity_extreme}
\end{figure}

Figure~\ref{fig:cavity_extreme} shows the temperature distribution for two limiting cases in the $(\kappa,\beta)$ parameter space, where $N=\kappa/\beta$ characterises the conduction-radiation balance. In the conduction-dominated limit ($\kappa=100$, $\beta=1$, $N=100$), the field resembles the pure-conduction solution: the temperature drops most steeply adjacent to the hot wall and then varies more gradually, with isotherms bending towards the cold side boundaries. In the radiation-dominated limit ($\kappa=0.01$, $\beta=100$, $N=10^{-4}$), in contrast, radiative energy emitted from the hot wall heats the medium throughout the cavity, producing a flatter, low-gradient temperature field with widely spaced isotherms. Both limiting behaviours agree with the profiles reported by Talukdar et al.~\cite{Talukdar2008} for equivalent values of $N$.

\begin{figure}[h]
    \centering
    \includegraphics[width=.45\textwidth]{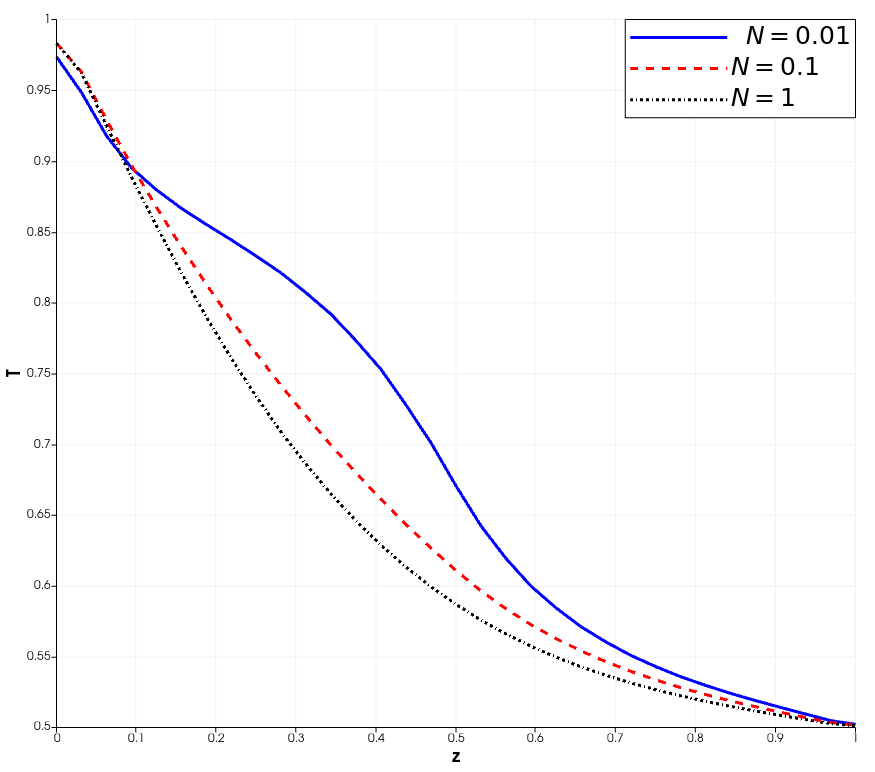}\hfill
    \includegraphics[width=.45\textwidth]{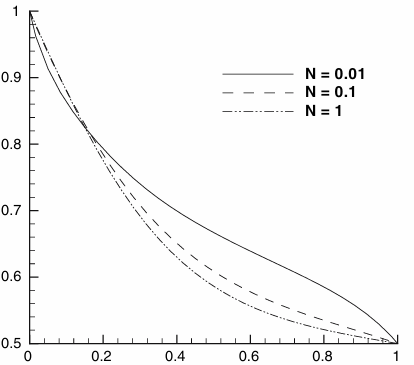}
    \caption{Centreline temperature profiles $(0.5,0.5,z)$ for $\beta=1$ and $\kappa=0.01,0.1,1.0$(corresponding to $N=\kappa/\beta=0.01,0.1,1$). Left: present computation. Right: reference profiles of Talukdar et al.~\cite{Talukdar2008} (Fig.~5a) at the same $N$ values.}
    \label{fig:cavity_Tz}
\end{figure}

Figure~\ref{fig:cavity_Tz} compares the centreline temperature profiles $(0.5,0.5,z)$ obtained with the present method for $\beta=1$ and $\kappa=0.01,0.1,1.0$ (corresponding to $N=\kappa/\beta=0.01,0.1,1$) with the reference profiles of Talukdar et al.~\cite{Talukdar2008} at the same $N$ values. In the radiation-dominated regime ($N=0.01$), volumetric radiative heating elevates the interior temperature, producing a relatively flat, low-gradient profile that lies above the pure-conduction limit. As $N$ increases, the profiles progressively approach the conduction-dominated solution, with $N=1$ nearly coinciding with the pure-conduction profile. The overall agreement with the reference data confirms the validity of the present numerical scheme.

\section{Conclusion}

This work proposes a monolithic solve and a block preconditioning method for the coupled nonlinear heat radiation
transport system.
The main principles are based on the approximation the physical coupling between temperature and radiation. $P_{\mathrm{schur}}$ retains full coupling and is the most robust across all regimes. $P_{\mathrm{split}}$ drops diffusion coupling and is tailored for the mixed regime. $P_{\mathrm{BJ}}$ neglects temperature-radiation coupling and is only efficient when heat conduction dominates.

A key insight from this study is that the performance of $P_{\mathrm{schur}}$ in the mixed regime is limited by the diagonal approximation of $L_T$, not by the Schur-complement framework itself. Exact Schur complement validation confirms this: the exact version restores the theoretical optimality of $P_{\mathrm{schur}}$ and performs comparably to $P_{\mathrm{split}}$. The split preconditioner thus offers a practical remedy by circumventing the diagonal approximation error entirely.

Numerical experiments on 3D problems confirm the hierarchy and demonstrate mesh-independent convergence across all regimes. Future work includes improved sparse approximations for the approximate Schur complements of radiation block, and extensions to higher-order discretisations and time-dependent problems.

\appendix
\section{Derivation of the dimensionless governing equations}
\label{sec:appendices}
In this appendix, quantities with an overbar  $\bar{ \ \ }$  denote dimensional quantities. The transient grey radiation transport equations in dimensional form are
\begin{align}
\frac{1}{c}\frac{\partial \bar{\Phi}}{\partial \bar{t}} + \Omega\cdot\nabla_{\bar{x}} \bar{\Phi} + \bar{\sigma}_a \bar{\Phi} - \frac{\bar{\sigma}_a}{\pi}\sigma_B \bar{T}^4 &= \bar{f},\quad \text{in }  D\times[0,T] \times\mathbb{S}^{2}, \label{A1} \\
\rho C_v \frac{\partial \bar{T}}{\partial \bar{t}} - \nabla_{\bar{x}}\cdot(\bar{\kappa} \nabla_{\bar{x}} \bar{T}) - \bar{\sigma}_a\oint_{\mathbb{S}^2} \bar{\Phi}\,\mathrm{d}\Omega + 4\bar{\sigma}_a\sigma_B \bar{T}^4 &= 0,\quad \text{in }  D\times[0,T]. \label{A2}
\end{align}

To non-dimensionalise the system, we introduce the characteristic scales
\begin{align}
&\Phi = \frac{\pi \bar{\Phi}}{\sigma_B T_{\mathrm{ref}}^4},\quad
T = \frac{\bar{T}}{T_{\mathrm{ref}}},\quad
\mathbf{x} = \frac{\bar{\mathbf{x}}}{L},\quad
t = \frac{\bar{t}}{t_{\mathrm{ref}}},\\
&\sigma = \bar{\sigma}_a L,\quad
\kappa = \frac{\bar{\kappa}}{\kappa_{\mathrm{ref}}},\quad
\beta = \frac{4\sigma_B T_{\mathrm{ref}}^3}{\kappa_{\mathrm{ref}} / L},\quad
f = \frac{\pi L \bar{f}}{\sigma_B T_{\mathrm{ref}}^4}. \label{A3}
\end{align}
Here $L$ is a characteristic length scale, $T_{\mathrm{ref}}$ is a reference temperature, and $\kappa_{\mathrm{ref}}$ is a reference thermal conductivity. Substituting \eqref{A3} into \eqref{A2} and dividing by $\kappa_{\mathrm{ref}} T_{\mathrm{ref}} / L^2$ gives
\[
\frac{\rho C_v L^2}{\kappa_{\mathrm{ref}} t_{\mathrm{ref}}}\frac{\partial T}{\partial t}
- \nabla_{\mathbf{x}}\cdot(\kappa \nabla_{\mathbf{x}} T)
- \beta\sigma\left(\frac{1}{4\pi}\oint_{\mathbb{S}^2}\Phi\,\mathrm{d}\Omega - T^4\right) = 0.
\]
Choosing $t_{\mathrm{ref}} = \rho C_v L^2 / \kappa_{\mathrm{ref}}$ makes the coefficient of the time-derivative term unity. Similarly, employing the same scaling to \eqref{A1} yields
\[
\frac{L}{c t_{\mathrm{ref}}}\frac{\partial \Phi}{\partial t}
+ \Omega\cdot\nabla_{\mathbf{x}}\Phi + \sigma(\Phi - T^4) = f.
\]

For the steady-state problems considered in this paper, all time-derivative terms are neglected, and the dimensionless steady system reads
\begin{align}
\Omega\cdot\nabla \Phi + \sigma\Phi - \sigma T^4 &= f,\quad \text{in } D\times \mathbb{S}^{2}, \label{A4} \\
-\nabla\cdot(\kappa \nabla T) - \beta\sigma\left(\frac{1}{4\pi}\oint_{\mathbb{S}^2}\Phi\,\mathrm{d}\Omega - T^4\right) &= 0,\quad \text{in } D. \label{A5}
\end{align}
Boundary conditions for $\Phi$ on inflow directions and for $T$ on $\partial D$ are imposed as in the main text.

\end{document}